\documentclass[11pt]{amsart}
\usepackage{hyperref}

\usepackage{enumerate}

\usepackage[utf8]{inputenc}
\usepackage[T1]{fontenc}

\usepackage[mathscr]{eucal}

\usepackage[a4paper, twoside=false, vmargin={2cm,3cm}, includehead]{geometry}

\newtheorem{lemma}{Lemma}[section]
\newtheorem{theorem}[lemma]{Theorem}
\newtheorem*{theorem*}{Theorem}

\newtheorem{proposition}[lemma]{Proposition}
\newtheorem*{proposition*}{Proposition}

\newtheorem{conjecture}{Conjecture}

\theoremstyle{definition}
\newtheorem*{claim*}{Claim}

\newtheorem*{definition}{Definition}

\newtheorem*{remark}{Remark}
\newtheorem*{remarks}{Remarks}

\newcommand{\C}{{\mathbb C}}

\newcommand{\N}{{\mathbb N}}

\newcommand{\R}{{\mathbb R}}

\newcommand{\Z}{{\mathbb Z}}

\newcommand{\one}{\mathbf{1}}

\newcommand{\norm}[1]{\left\Vert #1\right\Vert}

\begin{document}

\title[Multiple correlation sequences and nilsequences]{Multiple correlation sequences
and nilsequences}

\author{Nikos Frantzikinakis}
\address[Nikos Frantzikinakis]{University of Crete, Department of mathematics, Voutes University Campus, Heraklion 71003, Greece} \email{frantzikinakis@gmail.com}

\begin{abstract}
We study the structure of multiple correlation sequences defined by measure preserving actions of  commuting  transformations. When the iterates of the transformations are integer polynomials we prove  that any such correlation sequence is the sum of a nilsequence and an error term that is small in uniform density; this was previously known only for measure preserving actions of a single transformation.
We then use this decomposition result to give convergence criteria for multiple ergodic averages and deduce
some rather surprising results, for instance we infer
 convergence for actions of commuting transformations   from
the special case of actions of
 a single transformation. Our proof  of the decomposition result differs from previous works of V.~Bergelson,  B.~Host, B.~Kra, and A.~Leibman, as it does  not rely on the theory of characteristic factors.
  It consists of a  simple orthogonality argument and  the  main tool is an inverse theorem of B.~Host and B.~Kra for general bounded sequences.
\end{abstract}


\subjclass[2010]{Primary: 37A30; Secondary: 05D10,  11B30, 37A05. }

\keywords{Correlation sequences, nilsequences, multiple ergodic averages.}


\maketitle


\section{Introduction and results}
\subsection{Main result} Throughout this article, a \emph{system} is a probability space $(X,\mathcal{X},\mu)$  together with invertible, measure preserving
 transformations $T_1,\ldots, T_\ell\colon X\to X$  that commute. A multiple correlation sequence is a sequence of the form
$$
\int  T_1^{n_1}f_1\cdot \ldots \cdot T_\ell^{n_\ell}f_\ell\, d\mu
$$
where  $(X,\mathcal{X},\mu, T_1,\ldots, T_\ell)$ is a system,  $f_1,\ldots, f_\ell\in L^\infty(\mu)$, and $n_1,\ldots, n_\ell\in \Z$. The study of the limiting behavior of averages of such sequences, where the iterates are restricted to certain subsets of $\Z^\ell$,  has been an indispensable tool
in ergodic Ramsey theory and in particular in proving various far reaching extensions of Szemer\'edi's theorem on arithmetic progressions.
 Although the precise structure
  of the multiple correlation sequences
  is unknown even when $n_1=\cdots=n_\ell=n$, there is a widespread belief that modulo negligible terms the building blocks are sequences with algebraic structure (see
  \cite[Problem 1]{Fr11} for a related conjecture).
 \begin{definition}[\cite{BHK05}]
For $\ell\in \N$,  an  \emph{$\ell$-step nilsequence} is a sequence of the form $(F(g^n\Gamma))$, where $F\in C(X)$, $X=G/\Gamma$, $G$ is an  $\ell$-step nilpotent Lie group, $\Gamma$ is a discrete cocompact subgroup, and
  $g\in G$. A  \emph{$0$-step nilsequence}
is a constant sequence.
 \end{definition}
 When $T_i=T^i$, $i=1,\ldots, \ell,$ following the discovery of characteristic factors with  algebraic structure for some closely related multiple ergodic averages,
V.~Bergelson, B.~Host, and B.~Kra proved the following beautiful result (see also \cite{Mei90} for related work for $\ell=3$):
\begin{theorem*}[\mbox{\cite[Theorem 1.9]{BHK05}}]\label{T:0}
 For $\ell\in \N$, let $(X,\mathcal{X},\mu, T)$ be an ergodic system  and    $f_1,\ldots,f_\ell\in L^{\infty}(\mu)$ be functions with $\norm{f_i}_\infty\leq 1$. Then  we have the
 decomposition
$$
\int  T^{n}f_1\cdot \ldots \cdot T^{\ell n}f_\ell\ d\mu=a_{st}(n)+a_{er}(n), \quad  n\in \N, $$
where
\begin{enumerate}
\item $(a_{st}(n))$ is a uniform limit of  $(\ell-1)$-step nilsequences
with $\norm{a_{st}}_{\infty}\leq 1$;

\item  $\lim_{N-M\to \infty} \frac{1}{N-M}\sum_{n=M}^{N-1} |a_{er}(n)|^2=0$.
\end{enumerate}
\end{theorem*}
This result was  extended by A.~Leibman to cover  polynomial iterates in \cite{L10} and  not necessarily ergodic transformations in \cite{L14}.
The proofs of these results depend in an essential way on the fact that  characteristic factors for some suitable multiple ergodic averages
are inverse limits of nilsystems. This is no longer true for correlation sequences involving actions of  commuting transformations, which is  why efforts to prove decomposition results for such sequences  did not bring any results so far.
 In fact,  characteristic factors for commuting actions
are known to be extremely complex (for related work  see \cite{Au11a, Au11b}) which  has raised suspicions that  decomposition   results in this more general setup
   may  involve sequences very different
 from nilsequences. Our main result settles this rather elusive problem;
 we show that
modulo error terms that are small in uniform density,   correlation sequences of actions of commuting transformations are nilsequences.

\begin{theorem}\label{T:1}
 For $\ell\in \N$ let $(X,\mathcal{X},\mu, T_1,\ldots, T_\ell)$ be a system    and $f_1,\ldots,f_\ell\in L^{\infty}(\mu)$ be functions with $\norm{f_i}_\infty\leq 1$. Then for every $\varepsilon>0$ we have the
 decomposition
\begin{equation}\label{E:mulcom}
\int  T_1^{n}f_1\cdot \ldots \cdot T_\ell^{n}f_\ell\ d\mu=a_{st}(n)+a_{er}(n), \quad  n\in \N,
\end{equation}
where
\begin{enumerate}
\item $(a_{st}(n))$ is an $(\ell-1)$-step nilsequence with $\norm{a_{st}}_{\infty}\leq 1$;

\item \label{E:2} $\lim_{N-M\to \infty} \frac{1}{N-M}\sum_{n=M}^{N-1} |a_{er}(n)|^2\leq \varepsilon$.
\end{enumerate}
\end{theorem}
\begin{remark}
We do not know if a strengthening similar to the one in \cite[Theorem 1.9]{BHK05} holds where one
uses uniform limits of nilsequences in $(i)$ and
 takes $\varepsilon=0$ in $(ii)$.
\end{remark}
Our argument is rather versatile and does not rely on the theory of characteristic factors; we  rather focus on some distinctive properties correlation sequences as  in \eqref{E:mulcom} satisfy (see Theorem~\ref{T:2}). The idea that starts the proof comes from answering the following natural question:
``Can a multiple correlation sequence  as in \eqref{E:mulcom} be   asymptotically orthogonal to all $(\ell-1)$-step nilsequences?''.
On the one hand,  using an inverse theorem of B.~Host and B.~Kra
(see Theorem~\ref{T:inverse}), one gets that any such  sequence  has to be $U_{\ell}$-uniform.
On the other hand,  by successively applying van der Corput's lemma one sees that
a sequence of the form  \eqref{E:mulcom} is asymptotically orthogonal to all $U_{\ell}$-uniform sequences. Hence, any sequence that provides a positive answer to our question has to be asymptotically orthogonal to itself, that is,  has to converge to $0$ in density.

With this idea in mind, we  prove our main result as follows:  Given a sequence $(a(n))$ as in \eqref{E:mulcom}, we consider the
$(\ell-1)$-step nilsequence, call it $a_{st}$, that lies ``closest'' to $(a(n))$  with respect to the semi-norm $\norm{\cdot}_2$ defined  in \eqref{E:seminorm}. Then $a_{er}:=a-a_{st}$
is  asymptotically orthogonal to all $(\ell-1)$-step nilsequences, and arguing as before, we
get that $a_{st}$ and $a_{er}$ have the asserted properties.
A slight complication appears because  for $\ell\geq 2$  the  space of
$(\ell-1)$-step nilsequences (or uniform limits of such sequences)  is not $\norm{\cdot}_2$-complete; this is the reason why we are led to  an error term  $a_{er}$ that is small, but not zero,  in uniform density. For our argument to work we also have to make sure that various limits of uniform Ces\`aro
averages exist; to guarantee this, we use  a   result of T.~Austin~\cite{Au09}.

Using a variant of the previous argument  and a result of M.~Walsh~\cite{W12}
we get:

\begin{theorem}\label{T:1'}
 Let   $\ell,m\in \N$ and  $p_{i,j}\in \Z[t]$, $i=1,\ldots, \ell, j=1,\ldots, m$, be polynomials.
Then  there exists $k\in \N$, $k=k(\ell,m,\max{\deg(p_{i,j})})$,
such that for every system
  $(X,\mathcal{X},\mu, T_1,\ldots, T_\ell)$, functions
 $f_1,\ldots,f_m\in L^{\infty}(\mu)$ with $\norm{f_i}_\infty\leq 1$, and $\varepsilon>0$, we have
\begin{equation}\label{E:complicated}
\int  (\prod_{i=1}^\ell T_i^{p_{i,1}(n)})
f_1\cdot \ldots \cdot (\prod_{i=1}^\ell T_i^{p_{i,m}(n)})
f_m\, d\mu =a_{st}(n)+a_{er}(n), \quad n\in \N,
\end{equation}
where
\begin{enumerate}
\item $(a_{st}(n))$ is a $k$-step nilsequence with $\norm{a_{st}}_{\infty}\leq 1$;

\item \label{E:2} $\lim_{N-M\to \infty} \frac{1}{N-M}\sum_{n=M}^{N-1} |a_{er}(n)|^2\leq \varepsilon$.
\end{enumerate}
\end{theorem}

\subsection{A more general framework}
It turns out that Theorem~\ref{T:1} is a manifestation of a more general principle which
  asserts that if a sequence  is asymptotically orthogonal to all $U_{\ell}$-uniform sequences and
satisfies some necessary regularity conditions, then it  admits a decomposition like the one in Theorem~\ref{T:1}. To make this more precise we introduce  some notation (see  Section~\ref{SS:seminorms} for the definition of the uniformity seminorms).
\begin{definition}
Let $\ell\in \N$. We say that the bounded sequence $a\colon \N\to \C$ is
\begin{enumerate}
\item
 {\em $\ell$-anti-uniform}  if
 there exists
 $C:=C(\ell,a)$ such that
 $$
 \limsup_{N-M\to
  \infty} \Big|\frac{1}{N-M}\sum_{n=M}^{N-1} a(n)b(n)\Big|\leq C \norm{b}_{U_{\ell}(\mathbb{N})}$$ for every $b\in \ell^{\infty}$.

\item  {\em $\ell$-regular} if
 the limit
 $$
 \lim_{N-M\to \infty} \frac{1}{N-M}\sum_{n=M}^{N-1} a(n)\psi(n)
 $$ exists for every $(\ell-1)$-step nilsequence $(\psi(n))$.
\end{enumerate}
\end{definition}

 \begin{theorem}\label{T:2}
For $\ell\in \N$ let $a\colon \N\to \C$ be a sequence with $\norm{a}_\infty\leq 1$ that is $\ell$-anti-uniform and $\ell$-regular.
Then for every $\varepsilon>0$ we have the decomposition
$$
a(n)=a_{st}(n)+a_{er}(n), \quad n\in \N,
$$
where
\begin{enumerate}
\item \label{E:1a} $(a_{st}(n))$ is an $(\ell-1)$-step nilsequence with $\norm{a_{st}}_{\infty}\leq 1$;

\item \label{E:1b} $\lim_{N-M\to \infty} \frac{1}{N-M}\sum_{n=M}^{N-1} |a_{er}(n)|^2\leq \varepsilon$.
\end{enumerate}
\end{theorem}
\begin{remarks}
For general $\ell$-regular sequences a similar result is proved in \cite[Theorem 2.19]{HK09} with an error term that is small with respect to the seminorm  $\norm{\cdot}_{U_\ell(\N)}$.

A sequence $(a(n))$ that satisfies the asserted decomposition has to be $\ell$-regular.   It also
has to satisfy the estimate defining the $\ell$-anti-uniformity property if one introduces an arbitrarily small error term $\varepsilon$ on the right hand side and allows $C$ to  depend on $\varepsilon$ (this follows from
 \cite[Theorem 2.14]{HK09}).

Theorem~\ref{T:2} fails if we use standard Ces\`aro averages to define the notions of anti-uniformity and regularity (and leave the definition of $\norm{\cdot}_{U_\ell(\N)}$ as is); the sequence
$(e^{ i \sqrt{n}})$, illustrates this. The same sequence shows that
 anti-uniformity does not imply regularity ($(e^{ i \sqrt{n}})$ is $2$-anti-uniform but not $1$-regular).
\end{remarks}
\subsection{Applications}\label{SS:Applications}
On $\ell^\infty(\N)$ we define the seminorm $\norm{\cdot}_2$ by
\begin{equation}\label{E:seminorm}
\norm{a}_2^2:=\limsup_{N-M\to \infty} \frac{1}{N-M}\sum_{n=M}^{N-1}|a(n)|^2.
\end{equation}
 For $\ell\in \N$ we consider the following  subspaces of $\ell^\infty(\N)$:
$$
\mathcal{A}_\ell:=\Big\{(\psi(n))  \colon \psi \text{ is an } (\ell-1)\text{-step nilsequence}\Big\};
$$
$$
\mathcal{B}_\ell:=\Big\{\Big(\int   T^{k_1n}f_1\cdot \ldots \cdot T^{k_{\ell} n}f_\ell \, d\mu\Big)\colon  (X,\mathcal{X}, \mu,T) \text{ is a system}, f_i\in L^\infty(\mu),  k_i=\ell!/i \Big\};
$$
$$
\mathcal{C}_\ell:=\Big\{\Big(\int T_1^nf_1\cdot \ldots \cdot T_{\ell}^{n}f_\ell \, d\mu\Big)\colon (X,\mathcal{X}, \mu, T_1,\ldots, T_\ell) \text{ is a system and } f_i\in L^\infty(\mu)\Big\}.
$$
After Proposition~\ref{P:nilkey} we explain  why in the definition of $\mathcal{B}_\ell$ we  use  the  exponents $k_1,\ldots, k_\ell$ instead of  $1,\ldots, \ell$.
The space  $\mathcal{A}_\ell$ is linear since if  for $i=1,2$, $(F_i(g_i^n\Gamma_i))$  are $(\ell-1)$-step nilsequences on $G_i/\Gamma_i$, then their sum is  the $(\ell-1)$-step nilsequence $(F(g^n\Gamma))$ on $G/\Gamma$, where $G=G_1\times G_2$,
$\Gamma:=\Gamma_1\times \Gamma_2$, $g:=(g_1,g_2)$, $F(g\Gamma):=F_1(g_1\Gamma_1)+F_2(g_2\Gamma_2)$.
To see that the space
$\mathcal{C}_\ell$ is linear (similarly for $\mathcal{B}_\ell$), let $a,b\in \mathcal{C}_\ell$  be  defined by the systems $(X_i, \mathcal{X}_i,  \mu_i, T_i)$ and the functions $f^i_1,\ldots, f^i_\ell$, $i=1,2$. Then $c:=(a+b)/2$  is also a multiple correlation sequence  defined by the system $(X, \mathcal{X},  \mu, T)$, where $X=X_1\cup X_2$ (considered as disjoint subsets) with the corresponding $\sigma$-algebra $\mathcal{X}$, $\mu:=(\mu_1+\mu_2)/2$, $T$ equals $T_1$ on $X_1$ and $T_2$ on $X_2$, and $f_i:=f^1_i{\bf 1}_{X_1}+f^2_i{\bf 1}_{X_2}$, $i=1,\ldots, \ell$.

It is a rather striking fact that, modulo sequences that are small in uniform density,  the three subspaces
$\mathcal{A}_\ell$, $\mathcal{B}_\ell$,  $\mathcal{C}_\ell$
 coincide.
\begin{theorem} \label{T:3}
For every $\ell\in \N$ we have
$$
\overline{\mathcal{A}_{\ell}}=\overline{\mathcal{B}_\ell}=\overline{\mathcal{C}_\ell}
$$
where the closure is taken with respect to the seminorm $\norm{\cdot}_2$ defined in \eqref{E:seminorm}.
\end{theorem}
It is not hard to see that the first equality fails if we consider  closures with respect to the $\norm{\cdot}_\infty$ norm. The second equality may still hold under such circumstances but this is not something we can prove with the methods developed so far.

 The next two results  illustrate some rather surprising principles:
 $(i)$ convergence results for  actions of a single transformation automatically imply stronger convergence results for actions of commuting transformations; and
 $(ii)$ convergence results involving linear iterates
 automatically imply stronger convergence results involving polynomial iterates.
 \begin{theorem} \label{T:4}
Let $(r_n)$ be a strictly increasing  sequence of integers such that $r_n=O(n)$. Then for every $\ell\in \N$ the following statements  are equivalent:
\begin{enumerate}
\item \label{E:41} For every $(\ell-1)$-step nilsequence $(\psi(n))$ the limit $\lim_{N\to \infty}\frac{1}{N}\sum_{n=1}^{N}\psi(r_n)$ exists.

\item \label{E:42} For every system
$(X,\mathcal{X},\mu, T)$, functions $f_1,\ldots, f_\ell\in L^\infty(\mu)$, and for $k_i=\ell!/i$, $i=1,\ldots, \ell$,
 the following limit exists
 $$
 \lim_{N\to \infty}\frac{1}{N}\sum_{n=1}^{N} \int   T^{k_1r_n}f_1 \cdot\ldots\cdot   T^{k_\ell r_n}f_\ell \, d\mu.
 $$

\item  \label{E:43}For every system
$(X,\mathcal{X},\mu, T_1,\ldots, T_\ell)$ and functions $f_1,\ldots, f_\ell\in L^\infty(\mu)$
 the following limit exists
 $$
 \lim_{N\to \infty} \frac{1}{N}\sum_{n=1}^{N} \int   T_1^{r_n}f_1\cdot \ldots\cdot  T_{\ell}^{r_n} f_\ell \, d\mu.
 $$
\end{enumerate}
\end{theorem}
\begin{remark}
Equivalently, the growth condition $r_n=O(n)$ holds  if  the set $R:=\{r_1,r_2,\ldots\}$ has positive lower natural density.
\end{remark}
In the previous result we have established an equivalence for every fixed $\ell\in \N$, in the next result
we have to assume that a certain property is known for every $\ell\in \N$ in order to establish an equivalence (this is needed  for the equivalence of  (ii) and (iii)).
\begin{theorem} \label{T:4'}
Let $(r_n)$ be a strictly increasing  sequence of integers such that $r_n=O(n)$. Then  the following statements  are equivalent:
\begin{enumerate}
\item \label{E:4'1} For every
$\ell\in \N$ and $\ell$-step
nilsequence $(\psi(n))$ the limit $\lim_{N\to \infty}\frac{1}{N}\sum_{n=1}^{N}\psi(r_n)$ exists.

\item \label{E:4'2} For every  $\ell\in \N$, system
$(X,\mathcal{X},\mu, T)$, and functions $f_1,\ldots, f_\ell\in L^\infty(\mu)$,
 the following limit exists
 $$
 \lim_{N\to \infty}\frac{1}{N}\sum_{n=1}^{N} \int   T^{r_n}f_1\cdot \ldots \cdot T^{\ell r_n}f_\ell \, d\mu.
 $$

\item \label{E:4'3} For every $\ell\in \N$, polynomials $p_1,\ldots, p_\ell\in \Z[t]$,  system
$(X,\mathcal{X},\mu, T_1,\ldots, T_\ell)$, and functions $f_1,\ldots, f_\ell\in L^\infty(\mu)$,
 the following limit exists
 $$
 \lim_{N\to \infty}\frac{1}{N}\sum_{n=1}^{N} \int   T_1^{p_1(r_n)}f_1\cdot \ldots\cdot  T_{\ell}^{p_\ell(r_n)} f_\ell \, d\mu.
 $$
\end{enumerate}
\end{theorem}
Similar results hold if in  \eqref{E:4'1}-\eqref{E:4'3} of Theorems~\ref{T:4} and \ref{T:4'}   one replaces the limit
$\lim_{N\to \infty}\frac{1}{N}\sum_{n=1}^{N}$ with the limit $\lim_{N-M\to \infty}\frac{1}{N-M}\sum_{n=M}^{N-1}$  and the growth assumption on $(r_n)$ with the assumption that the range of this sequence has positive lower Banach density. Furthermore, the same method can be used to prove   convergence criteria  for weighted  averages where for a given
 bounded sequence of complex numbers $(w_n)$ one 
replaces  in \eqref{E:4'1}-\eqref{E:4'3} of Theorems~\ref{T:4} and \ref{T:4'}  the   averaging operation  $\frac{1}{N}\sum_{n=1}^{N}$
with the  averaging operation $\frac{1}{N}\sum_{n=1}^{N} w_n$.

\subsection{Conjectures}
 The growth assumption on $(r_n)$  in Theorems~\ref{T:4} and  \ref{T:4'}  is crucial for our argument to work as the proofs  use
Theorem~\ref{T:1} which is not helpful  for   sequences that grow faster than linearly. Nevertheless, we
 believe  that the following is true:
\begin{conjecture}
In Theorems~\ref{T:4} and  \ref{T:4'}  the growth assumption on  $(r_n)$  is superfluous.
\end{conjecture}
We also believe in the following strengthening of the second identity in Theorem~\ref{T:3}:
\begin{conjecture}
For every $\ell\in \mathbb{N}$ we have
$\overline{\mathcal{B}_\ell}=\overline{\mathcal{C}_\ell}$ where the closure is taken with respect to the norm $\norm{\cdot}_\infty$.
\end{conjecture}

\subsection{Notation}
  We denote by $\N$ the set of positive integers.

\noindent
If $(a(n))$ is a bounded sequence we denote by  $\limsup_{N-M\to \infty} |\frac{1}{N-M}\sum_{n=M}^{N-1}a(n)|$ the limit (it exists by subadditivity)
$
\lim_{N\to \infty} \sup_{M\in \N} \Big|\frac{1}{N}\sum_{n=M}^{M+N-1}a(n)\Big|.
$


\subsection{ Acknowledgements.} I would like to thank B.~Host, B.~Kra,  M.~Wierdl, and the referee for helpful remarks.

\section{Proofs of  results}
\subsection{Uniformity seminorms and the Host-Kra inverse theorem}\label{SS:seminorms}
We give  a slight variant of the uniformity seminorms defined by B.~Host and B.~Kra in \cite{HK09}.
\begin{definition}
Let $\ell\in \N$ and $a\colon \N\to \C$ be a bounded sequence.

\begin{enumerate}
\item
Given a sequence of intervals  ${\bf I}=(I_N)$   with  lengths tending to infinity,
we say that the   sequence $(a(n))$ is {\em distributed
 regularly along ${\bf I}$}  if the limit
 $$
 \lim_{N\to \infty} \frac{1}{|I_N|}\sum_{n\in I_N} a_1(n+h_1)\cdot\ldots\cdot a_r(n+h_r)
$$
exists for every $r\in \N$ and  $h_1,\ldots, h_r\in \N$, where $a_i$ is either $a$ or $\bar{a}$.

\item If ${\bf I}$  is as in (i)  and $(a(n))$ is  distributed
 regularly along ${\bf I}, $ we define inductively
$$\norm{a}_{{\bf I}, 1}:=
\lim_{N\to \infty} \Big|\frac{1}{|I_N|}\sum_{n\in I_N} a(n)\Big|;
$$
and for $\ell\geq 2$ (one  can  show as in \mbox{\cite[Proposition 4.3]{HK09}} that the next limit exists)
$$
\norm{a}_{{\bf I}, \ell}^{2^{\ell}} :=\lim_{H\to \infty} \frac{1}{H}\sum_{h=1}^H \norm{\sigma_ha\cdot \bar{a}}^{2^\ell-1}_{{\bf I}, \ell-1}
$$
where $\sigma_h$ is the shift transformation defined by $(\sigma_ha)(n):=a(n+h)$.
\item If $(a(n))$ is a bounded sequence we let
$$
\norm{a}_{U_\ell(\N)}:=\sup_{{\bf I}}\norm{a}_{{\bf I}, \ell}
$$
where the sup is taken over all sequences of intervals ${\bf I}$  with lengths tending to infinity along which
the sequence  $(a(n))$ is  distributed
 regularly.
\end{enumerate}
\end{definition}
 An  application of Lemma~\ref{L:VDC}
shows that $\norm{a}_{{\bf I}, 1}$, as defined  here, is smaller than the corresponding quantity
defined in \cite{HK09} (they can be different though). Furthermore, the inductive formula is identical in both cases
(see \cite[Proposition 4.4]{HK09}), hence  $\norm{\cdot}_{U_\ell(\N)}$,  as defined here, is a seminorm that is  smaller than the corresponding seminorm
defined in \cite{HK09}. In fact, it can be shown that the two seminorms coincide but we will not need this.

 Using the main structural result in   \cite{HK05},  B.~Host and B.~Kra  proved an  inverse theorem
that will be  a key ingredient in  the proof of Theorem~\ref{T:2}. We state a slight variant of it next (\cite[Theorem 2.16]{HK09}
 gives a stronger lower bound but it does not allow to assume that $\norm{b}_\infty\leq 1$). Its proof
 amounts to a simple modification of the argument given in \cite[Theorem 2.16]{HK09}; we give the details for completeness.
\begin{theorem}[\mbox{\cite[Theorem 2.16]{HK09}}]\label{T:inverse}
Let $a\colon \N\to \C$ be a sequence  of complex numbers  with $\norm{a}_\infty\leq 1$ and $\ell\in \N$. Then for every $\varepsilon>0$ there exists an $(\ell-1)$-step  nilsequence $(b(n))$ with $\norm{b}_\infty\leq 1$
 such that
$$
\limsup_{N-M\to \infty}\Big|\frac{1}{N-M}\sum_{n=M}^{N-1} a(n)b(n)\Big|\geq \norm{a}_{U_{\ell}(\N)}^{2^\ell}-\varepsilon.
$$
\end{theorem}
\begin{remark}
It is crucial that the  seminorms were defined using uniform and not standard Ces\`aro averages as  in the latter case it  is shown in \mbox{\cite[Paragraph 2.4.3]{HK09}} that the corresponding inverse theorem fails.
For  standard Ces\`aro averages
a finitary inverse theorem was proved in~\cite{GTZ12c} but it is not clear whether it has   an infinitary variant that  is useful for our purposes.
\end{remark}
\begin{proof}
We refer the reader to \cite{HK09} for notation used in this argument. In what follows
we assume that the seminorms
 $\norm{a}_{{\bf I}, \ell}$ are defined as in \cite{HK09}.

Let $0<\varepsilon<1$. By \cite[Proposition 6.2]{HK09} there exists a sequence of intervals ${\bf I}=(I_N)$ with lengths tending to infinity and an  $(\ell-1)$-step nilsequence  $(c(n))$ of the form $c(n)=F(g^n\Gamma)$, where $F$ is a continuous function on an $(\ell-1)$-step nilmanifold
$X=G/\Gamma$ and $g\in G$ is an element that acts ergodically on $X$, such that the sequences $a-c$ and $a$ satisfy property $\mathcal{P}(\ell)$ on ${\bf I}$ and moreover   we have the estimates
 \begin{equation}\label{E:2est}
 \norm{a-c}_{{\bf I}, \ell}\leq \varepsilon, \quad \norm{a}_{{\bf I}, \ell}\geq \norm{a}_{U_\ell(\N)}-\varepsilon.
 \end{equation}
Furthermore,
 we have   $\norm{F}_\infty\leq 1$, this is because in the proof of  \cite[Proposition 6.2]{HK09} the function $F$ is defined as a conditional expectation of a function bounded by $1$.  We let $b(n):=H(g^n\Gamma)$, where $H:=\mathcal{D}_\ell F$, and check that the asserted properties are satisfied.

 First note that $(b(n))$ is an $(\ell-1)$-step nilsequence  and  since $\norm{F}_\infty\leq 1$ we have $\norm{H}_\infty\leq 1$, hence $\norm{b}_\infty\leq 1$. Furthermore, by  \cite[Corollary 5.3]{HK09} we have  $H\in C(X)$, hence $F\cdot H\in C(X)$,  and since  $g$ acts ergodically on $X$ we have
  $$
  \lim_{N\to \infty}\frac{1}{|I_N|}\sum_{n\in I_N} c(n)b(n)=
  \int F\cdot H\, dm_X=
  \norm{F}_\ell^{2^\ell}=\norm{c}_{{\bf I}, \ell}^{2^\ell}
  $$
  where we used  the identity $\int F\cdot \mathcal{D}_\ell F\, dm_X=\norm{F}_\ell^{2^\ell}$ and
  \cite[Corollary 3.11]{HK09} to justify the last two identities.
  By \eqref{E:2est} and the triangle inequality this is greater or equal than
  $$
  (\norm{a}_{{\bf I}, \ell}-\varepsilon)^\ell\geq
  (\norm{a}_{U_\ell(\N)}-2\varepsilon)^\ell
  \geq \norm{a}_{U_\ell(\N)}^\ell-k_\ell\varepsilon
  $$
  for some positive integer $k_\ell$.
  On the other hand, by \cite[Theorem 2.13]{HK09} we have
  $$
  \limsup_{N\to \infty}\Big|\frac{1}{|I_N|}\sum_{n\in I_N}(a(n)-c(n))b(n) \Big|\leq
  \norm{a-c}_{{\bf I}, \ell} \norm{b}^*_{\ell}\leq \varepsilon
  $$
  where we used \eqref{E:2est} and that  $\norm{b}^*_{\ell}=\norm{\mathcal{D}_\ell F}_\ell^*=\norm{F}_\ell^{2^\ell-1}\leq 1$
 (the second identity follows from  \cite[Equation (14)]{HK09}). Combining the previous bounds we get the asserted result.
\end{proof}

\subsection{Proof of Theorem~\ref{T:2}}
  Let $\ell\in \N$ and  $(a(n))$ be an $\ell$-regular and $\ell$-anti-uniform sequence with $\norm{a}_\infty\leq 1$.
 We first remark that the limit
 \begin{equation}\label{E:exists}
  \lim_{N-M\to \infty} \frac{1}{N-M}\sum_{n=M}^{N-1} |a(n)|^2 \quad \text{ exists}.
 \end{equation}
 This follows from our anti-uniformity assumption
 and  \cite[Theorem 2.19]{HK09} (it applies since $(a(n))$ is $\ell$-regular) which states
 that for every $\epsilon>0$ we have a decomposition $a=a_1+a_2$ where $a_1$ is an $(\ell-1)$-step nilsequence
 and $\norm{a_2}_{U_\ell(\N)}\leq \epsilon$. Writing $|a(n)|^2=a\bar{a}_1+a\bar{a}_2$ one  checks the asserted convergence at once.

  We let
 $$
Y:=\Big\{(\psi(n)) \colon \psi \text{ is an } (\ell-1)\text{-step nilsequence}\Big\}
$$
and
$$
X:=\text{span}\{{Y,a\}}.
$$
On $X\times X$ we define the bilinear form
$$
\langle f, g \rangle:=\lim_{N-M\to \infty} \frac{1}{N-M}\sum_{n=M}^{N-1} f(n)\overline{g}(n).
$$
  Note that the limit exists for $f,g\in X$. This is the case  if $f$ or $g$ is equal to $a$ because of  our regularity assumption and \eqref{E:exists}, and when both $f$ and $g$ are in $Y$ because  limits of
uniform Ces\`aro averages of nilsequences exist \cite{L05, Les91}. This bilinear form induces the  seminorm
$$
\norm{f}_2:=\sqrt{\langle f, f \rangle}.
$$
This is the restriction on $X$ of the seminorm  \eqref{E:seminorm} defined on $\ell^\infty(\N)$.

Let $\varepsilon>0$.
 There exists  $y_0\in Y$  such that
\begin{equation}\label{E:distance}
\norm{a-y_0}_2^2\leq d^2+\delta^2
\end{equation}
where
\begin{equation}\label{E:ded}
 d:=\inf\{ \norm{a-y}_2\colon y\in Y\}, \quad \delta:=(\varepsilon/(4C))^{2^\ell},
 \end{equation}
and  $C:=C(\ell,a)$ is the constant determined by our $\ell$-anti-uniformity assumption on $a$. We can assume that $C\geq 1$.
Furthermore,  we can assume without loss of generality that
\begin{equation}\label{E:y0}
\norm{y_0}_\infty\leq 1.
\end{equation}
Indeed, let $y_0:=(F(g^n\Gamma))$ where $X=G/\Gamma$ is a nilmanifold, $g\in G$, and $F\in C(X)$.
Then  the sequence $\tilde{y}_0:=(\tilde{F}(g^n\Gamma))$,  where $\tilde{F}:=F\cdot \one_{|F|\leq 1}+e^{2\pi i \arg(F)}\cdot \one_{|F|\geq 1} \in C(X)$, is a nilsequence, $\norm{\tilde{y}_0}_\infty\leq 1$,  and
as $\norm{a}_\infty\leq 1$ we get that $|a(n)-\tilde{y}_0(n)|\leq |a(n)-y_0(n)|$ for every $n\in \N$, hence $\norm{a-\tilde{y}_0}_2\leq \norm{a-y_0}_2$.

 It follows from \eqref{E:distance} that for every $y\in Y$ we have
$$
-\delta^2\leq \norm{a-(y_0+\delta y)}_2^2-\norm{a-y_0}_2^2=-2\delta\text{Re}(\langle a-y_0,y\rangle)+\delta^2\norm{y}_2^2.
$$
 Hence,
 $$
\text{Re}(\langle a-y_0,y\rangle) \leq \delta \ \text{ for every } y\in Y \text{ with } \norm{y}_2\leq 1.
 $$
Inserting   $-y$ and $\pm i y$ in place of $y$ we deduce that
 \begin{equation}\label{E:ay0}
\sup_{y \in Y\colon \norm{y}_2\leq 1} |\langle a-y_0,y\rangle| \leq 2\delta.
 \end{equation}
 Since the set $\{y\in Y\colon \norm{y}_2\leq 1\}$ contains all $(\ell-1)$-step nilsequences that are bounded by $1$, we deduce  from Theorem~\ref{T:inverse}
that
\begin{equation}\label{E:uniform}
\norm{a-y_0}_{U_{\ell}(\N)}\leq (2\delta)^{2^{-\ell}}.
\end{equation}

We let
$$
a_{st}:=y_0, \quad a_{er}:=a-y_0.
$$
Then
$$
a=a_{st}+a_{er}
$$
and $(a_{st}(n))$ is an $(\ell-1)$-step nilsequence with $\norm{a_{st}}_\infty\leq 1$ by \eqref{E:y0}.
Since $a$ is $\ell$-anti-uniform we get using \eqref{E:uniform} and the definition of $\delta$ in \eqref{E:ded} that
  $$
  |\langle a,a_{er}\rangle|\leq C\norm{a_{er}}_{U_{\ell}(\N)}\leq \varepsilon/2.
  $$
Furthermore,   \eqref{E:ay0} gives
  $$
  |\langle a_{st},a_{er}\rangle|\leq \varepsilon/2.
  $$
 Combining the last  two estimates  we deduce that
$$
 \norm{a_{er}}_2^2 =\langle a_{er},a_{er}\rangle\leq
 |\langle a,a_{er}\rangle|+ |\langle a_{st},a_{er}\rangle|
\leq \varepsilon.
$$
  This completes the proof of Theorem~\ref{T:2}.




  \subsection{Proof of Theorem~\ref{T:1}} In view of  Theorem~\ref{T:2}, it suffices to prove that for every $\ell\in \N$ the sequence
  $a\colon \N\to \C$ defined by
\begin{equation}\label{E:an}
a(n):=\int T_1^{n}f_1\cdot \ldots \cdot T_\ell^{n}f_\ell\ d\mu, \quad n\in \N,
\end{equation}
is $\ell$-anti-uniform and $\ell$-regular.

\subsubsection{Anti-uniformity} Throughout, we can and will assume that $\norm{f_i}_\infty\leq 1$ for $i=1,\ldots, \ell$.
The $\ell$-anti-uniformity   follows by successive applications of the following Hilbert space variant of van der Corput's estimate
(for a proof see \cite{Be87a}).
 \begin{lemma}\label{L:VDC}
Let  $(v_n)$ be a bounded  sequence of vectors in an inner product
space and $(I_N)$ be a sequence of intervals with lengths tending to infinity.
Then
$$
\limsup_{N\to\infty}
\norm{\frac{1}{|I_N|}\sum_{n\in I_N} v_n}^2\leq 4 \ \!
\limsup_{H\to\infty} \frac{1}{H}\sum_{h=1}^H
\limsup_{N\to\infty}\Big|
\frac{1}{|I_N|}\sum_{n\in I_N} \langle v_{n+h},v_{n}\rangle \Big|.
$$
\end{lemma}
It suffices  to show that    for every $\ell\in \N$ and  every sequence of intervals ${\bf I}:=(I_N)$ with lengths tending to infinity,
any sequence $(a(n))$ given by \eqref{E:an}  satisfies the estimate
$$
 \limsup_{N\to
  \infty} \Big|\frac{1}{|I_N|}\sum_{n\in I_N} a(n)b(n)\Big|\leq 4 \norm{b}_{U_{\ell}(\N)}
  $$
  for every $b\in \ell^{\infty}(\N)$. Using a diagonal argument and passing to a subsequence  of $(I_N)$
   (if necessary) we can and will assume that the sequence $(b(n))$ is distributed regularly along  the sequence
  ${\bf I}$.
  It suffices to establish that for any sequence $(a(n))$ as in \eqref{E:an} which is bounded by $1$  and any $b\in \ell^\infty(\N)$ which is distributed regularly along a sequence of intervals ${\bf I}$, we have
\begin{equation}\label{E:anbn}
 \limsup_{N\to
  \infty} \Big|\frac{1}{|I_N|}\sum_{n\in I_N} a(n)b(n)\Big|\leq 4  \norm{b}_{{\bf I}, \ell}.
  \end{equation}

We prove this  by induction on  $\ell$. For $\ell=1$ the result holds trivially.
Suppose that $\ell\geq 2$ and  the statement holds for $\ell-1$.
We compose with $T_\ell^{-n}$,  use the Cauchy-Schwarz inequality, and then  Lemma~\ref{L:VDC} (on the space $L^2(\mu)$) for the sequence
$$
v_n := b(n)\cdot \tilde{T}_1^nf_1\cdot \tilde{T}_2^nf_2\cdot\ldots\cdot \tilde{T}_{\ell-1}^nf_{\ell-1}, \quad n\in \N,
$$
where $\tilde{T}_i:=T_iT_\ell^{-1}$ for $i=1,\ldots,\ell-1$.
We deduce that the square of the left hand side in \eqref{E:anbn} is bounded by
\begin{equation}\label{E:u_n}
 \limsup_{N\to\infty}\Bigl\lVert \frac{1}{|I_N|}\sum_{n\in I_N}
v_{n}\Bigr\lVert_{L^2(\mu)}^2\leq 4 \limsup_{H\to
\infty}\frac{1}{H}\sum_{h=1}^{H} \limsup_{N\to\infty}\Big|
\frac{1}{|I_N|}\sum_{n\in I_N} \langle{ v_{n+h}, v_n \rangle  }\Big|.
\end{equation}
A simple computation gives that
$$
\frac{1}{|I_N|}\sum_{n\in I_N} \langle{ v_{n+h}, v_n \rangle }
=  \frac{1}{|I_N|}\sum_{n\in I_N} b(n+h)\cdot \bar{b}(n) \int
\tilde{T}_1^n \tilde{f}_{1,h} \cdot\ldots\cdot \tilde{T}_{\ell-1}^n
\tilde{f}_{\ell-1,h}\,d\mu
$$
where
$\tilde{f}_{j,h}=\tilde{T}_{j}^hf_{j}\cdot\bar{f}_{j}$ for $j=1,\ldots, \ell-1$.
 Note that the   maps
$\tilde{T}_1,  \ldots, \tilde{T}_{\ell-1}$ commute, for $h\in \N$ the sequence $(b(n+h) \bar{b}(n))$ is distributed regularly along ${\bf I}$,
 and $\norm{\tilde{f}_{j,h}}_\infty\leq 1$ for $j=1,\ldots, \ell-1$. Using the induction hypothesis and the defining property of the seminorms  we can bound the right hand
side in \eqref{E:u_n} by $16$ times
$$
\lim_{H\to \infty} \frac{1}{H}\sum_{h=1}^{H}\norm
{\sigma_hb\cdot b}_{{\bf I}, \ell-1 } \leq \lim_{H\to \infty}
\Big(\frac{1}{H} \sum_{h=1}^{H}\norm
{\sigma_hb\cdot b}_{{\bf I}, \ell-1 }^{2^{\ell-1}}\Big)^{1/2^{\ell-1}} =\norm {b}_{{\bf I}, \ell}^2
$$
where $(\sigma_h b)(n):=b(n+h)$.
Taking square roots we get the asserted estimate.

\subsubsection{Regularity} Let $\ell\in \N$.
To prove that $(a(n))$ is  $\ell$-regular   we will use
a  known mean convergence result for multiple ergodic averages and Proposition~\ref{P:nilkey} below.
 We start with the following
result of B.~Green and T.~Tao:
\begin{lemma}[\mbox{\cite[Lemma~14.2]{GT08}}]\label{L:nilkey}
For $\ell\in \N$ let $X=G/\Gamma$ be an $(\ell-1)$-step nilmanifold. Then there exists
a continuous map $P\colon X^{\ell}\to X$ such that
\begin{equation}\label{E:asd}
P(hg\Gamma, h^2g\Gamma, \ldots, h^{\ell} g\Gamma)=g\Gamma, \quad \text{ for every }  g,h\in G.
\end{equation}
\end{lemma}
The result in \mbox{\cite[Lemma~14.2]{GT08}} gives $P(g\Gamma,hg\Gamma, h^2g\Gamma, \ldots, h^{\ell-1} g\Gamma)=h^{\ell} g\Gamma$.
Inserting $h^{-\ell}g $ in place of $g$, then $h^{-1}$ in place of $h$, and rearranging coordinates,  we get  \eqref{E:asd}.
\begin{proposition}\label{P:nilkey}
For $\ell\in \N$ let  $(\psi(n))$ be an $(\ell-1)$-step nilsequence.
Then for every $\varepsilon>0$ there exists a
system $(X,\mathcal{X},\mu,T)$ and functions $f_1,\ldots, f_{\ell}\in L^\infty(\mu)$,
 such  that the
sequence $(b(n))$, defined by
\begin{equation}\label{E:bn}
b(n):=\int  T^{k_1n}f_1 \cdot \ldots \cdot  T^{k_{\ell} n}f_{\ell}\ d\mu, \quad n\in \N,
\end{equation}
where $k_i:=\ell!/i$ for $i=1,\ldots, \ell$, satisfies
$$
\norm{\psi-b}_\infty\leq \varepsilon.
$$
\end{proposition}
\begin{remarks}
To prove a variant of this result that uses the integers $1, \ldots, \ell$ in place of $k_1,\ldots, k_\ell$, 
one would have  to prove a non-trivial variant of Lemma~\ref{L:nilkey} that establishes in  place of \eqref{E:asd} the identity
$P(h^{k_1}g\Gamma, h^{k_2}g\Gamma, \ldots, h^{k_\ell} g\Gamma)=g\Gamma
$ for every $g,h\in G$.

Combining    \cite[Theorem A (ii)]{BL07} with Proposition~\ref{P:nilkey}  one deduces that for every bounded generalized polynomial $p\colon \N\to \R$ (see definition in \cite{BL07}) the sequences $(p(n))$
 and $(e^{ip(n)})$ can be approximated arbitrarily well  in $\norm{\cdot}_2$  by a sequence of the form \eqref{E:bn}.
\end{remarks}
\begin{proof}


Let $\varepsilon>0$ and
$$
\psi(n):=F(g^n\Gamma)
$$ where   $F\in C(X)$, $X=G/\Gamma$ is an $(\ell-1)$-step nilmanifold, and $g\in G$.

  By \cite[Paragraph 1.11]{L05}  we have that $X$ is isomorphic to a subnilmanifold of a nilmanifold $\tilde{X}=\tilde{G}/\tilde{\Gamma}$, where $\tilde{G}$ is a connected and simply connected $(\ell-1)$-step nilpotent Lie group, $\tilde{\Gamma}$ is a discrete cocompact subgroup of $\tilde{G}$, and all elements of $G$ are represented in $\tilde{G}$.  Then $\psi(n)=\tilde{F}(\tilde{b}^n\tilde{\Gamma})$ for some
$\tilde{b}\in \tilde{G}$ and $\tilde{F}\in C(\tilde{X})$. Hence, in what follows we can and will assume that the group $G$ is connected.

 Using Lemma~\ref{L:nilkey} with $g^n$ in place of $g$ and  $h:=g^m$, $m,n\in \N$,   we get that there exists a continuous map $P\colon X^{\ell}\to X$ such that
\begin{equation}\label{E:gn}
g^n\Gamma=P(g^{m+n}\Gamma, g^{2m+n}\Gamma,\ldots, g^{\ell m+n}\Gamma) \quad \text{ for every } m,n\in \N.
\end{equation}
Let $g_0\in G$ be  such that $g_0^{\ell!}=g$ (such a $g_0$ exists since $G$ is connected, hence  divisible) and for $i=1,\ldots, \ell$  let $g_i:=g_0^i$.
Applying \eqref{E:gn} with $g_0$ in place of $g$ and $\ell! n$ (a multiple of $n$ is needed that is divisible by all the coefficients of $m$ that appear in \eqref{E:gn}) in place of $n$ we get
$$
\psi(n)=F(g_0^{\ell! n}\Gamma)=\tilde{F}(g_1^{m+k_1n}\Gamma, g_2^{m+k_2n}\Gamma,\ldots, g_{\ell}^{ m+k_{\ell} n}\Gamma) \quad \text{ for every } m,n\in \N,
$$
where $\tilde{F}:=F\circ P\in C(X^{\ell})$.
Averaging over $m\in \N$  we get
 $$
\psi(n)=\lim_{M\to \infty} \frac{1}{M}\sum_{m=1}^M \tilde{F}(g_1^{m+k_1n}\Gamma, g_2^{m+k_2n}\Gamma,\ldots, g_{\ell}^{ m+k_{\ell} n}\Gamma) \quad \text{ for every }  n\in \N.
$$
Since $\tilde{F}$ can be approximated uniformly
by  linear combinations of  functions of the form $\tilde{f}_1\otimes\cdots\otimes \tilde{f}_{\ell}$, where for $i=1,\ldots, \ell$ the function $\tilde{f}_i \in C(X^\ell)$ depends on the coordinate $x_i$ only,  we get that $(\psi(n))$ can be  approximated in the $\norm{\cdot}_\infty$ norm within $\varepsilon$  by a finite linear combination of sequences $(a(n))$ of the form
\begin{equation}\label{E:limform}
a(n):=\lim_{M\to \infty} \frac{1}{M}\sum_{m=1}^M \tilde{f}_1(\tilde{g}^{m+k_1n}\tilde{\Gamma})\cdot  \tilde{f}_2(\tilde{g}^{m+k_2n}\tilde{\Gamma})\cdot \ldots \cdot \tilde{f}_{\ell}(\tilde{g}^{ m+k_{\ell} n}\tilde{\Gamma}), \quad n\in \N,
\end{equation}
where $\tilde{X}:=X^{\ell}$ , $\tilde{\Gamma}:=\Gamma\times\cdots \times \Gamma$, $\tilde{f}_i\in C(\tilde{X})$, and $\tilde{g}:=(g_1, \ldots, g_{\ell})$.
It is known (see \cite{L05} for example) that the limit in \eqref{E:limform} is equal to
$$
\int_{\tilde{Y}} \tilde{f}_1(\tilde{g}^{k_1n}\tilde{y})\cdot \tilde{f}_2(\tilde{g}^{k_2n}\tilde{y})\cdot \ldots \cdot \tilde{f}_{\ell}(\tilde{g}^{k_{\ell} n}\tilde{y}) \, dm_{\tilde{Y}}, \quad n\in \N,
$$
where $\tilde{Y}$ is the subnilmanifold of $\tilde{X}$ defined by  the closure of the set
$\{\tilde{g}^m \tilde{\Gamma}\colon  m\in \N\}$. This proves that  the sequence $(a(n))$   has the form \eqref{E:bn}.
  Since finite linear combinations of sequences of the form
\eqref{E:bn} still have the form \eqref{E:bn} (see Section~\ref{SS:Applications}) the  proof is complete.
\end{proof}

We are now ready to  verify that if  $(a(n))$ is as in  \eqref{E:an}, then it  is  $\ell$-regular for every $\ell\in \N$.
By Proposition~\ref{P:nilkey}, in order to check that the limit
 $\lim_{N-M\to \infty} \frac{1}{N-M}\sum_{n=M}^{N-1}a(n)\psi(n) $ exists for every $(\ell-1)$-step nilsequence $(\psi(n))$,
 it suffices to check that the limit
  \begin{equation}\label{E:abnn}
\lim_{N-M\to \infty} \frac{1}{N-M}\sum_{n=M}^{N-1} a(n)b(n)
\end{equation}
 exists for every sequence
  $(b(n))$ of the form $\int  S^{k_1n}g_1\cdot \ldots \cdot S^{k_{\ell} n}g_{\ell}\ d\nu$ , where $k_1,\ldots, k_\ell\in \N$,  $(Y,\mathcal{Y},\nu, S)$ is a system, and $g_1,\ldots, g_{\ell}\in L^\infty(\nu)$. This follows  from  the mean convergence result of T.~Austin~\cite{Au09} (which  strengthens  the convergence result of T.~Tao
   \cite{Ta08} to uniform averages)
  applied to
the transformations
 $\tilde{T}_i:=T_i\times S^{k_i}$ acting on $X\times Y$ with the measure $\tilde{\mu}:=\mu\times \nu$ and the functions $\tilde{f_i}:=f_i\otimes g_i \in L^\infty(\tilde{\mu})$, $i=1,\ldots, \ell$.

\subsection{Proof of Theorem~\ref{T:1'}}
Modulo a known convergence result of M.~Walsh~\cite{W12}
the argument is similar to the one used to prove Theorem~\ref{T:1}, we explain the minor modifications needed next.

To verify $k$-anti-uniformity for some $k\in \N$ that depends only on $\ell, m$ and the maximum degree of the polynomials $p_{i,j}$,
  one has to make successive uses  of Lemma~\ref{L:VDC}
and apply
an inductive argument, often called PET induction, introduced by V.~Bergelson  in \cite{Be87a}. The details are very similar  to those in the proof of \mbox{\cite[Lemma 3.5]{FrHK11}} and so we omit them.

To verify regularity, we can argue as in the case of linear iterates, using the convergence result of M.~Walsh  \cite{W12} for averages of expressions of the form \eqref{E:complicated}. At the very last step  one needs to verify that if $(a(n))$ is as in \eqref{E:complicated}, then the limit
\eqref{E:abnn}
 exists for every sequence
  $(b(n))$ of the form $\int  S^{k_1n}g_1\cdot \ldots \cdot S^{k_r n}g_r\ d\nu$,
   where $r\in \N$ is arbitrary, $k_1,\ldots, k_r\in \N$, $(Y,\mathcal{Y},\nu, S)$ is a system, and $g_1,\ldots, g_r\in L^\infty(\nu)$.
   The only change needed is to use Walsh's convergence result   for
the $\ell+r$ commuting measure preserving  transformations
$T_i\times \text{id}$, $i=1,\ldots, \ell$, and  $\text{id}\times S^{k_j}$, $j=1,\ldots, r$,  acting on $X\times Y$ with the measure $\tilde{\mu}:=\mu\times \nu$, and the functions $f_i\otimes 1$, $i=1,\ldots, \ell$  and $1\otimes g_j$, $j=1 \ldots, r$. If the polynomial iterates are chosen appropriately, one verifies  that  $a(n)b(n)$ is  also a multiple correlation sequence with polynomial iterates, hence, by Walsh's convergence result \cite{W12},  the limit \eqref{E:abnn} exists.

\subsection{Extension to nilpotent groups} Essentially the same argument can be used when the transformations $T_1,\ldots, T_\ell$ generate a nilpotent group; the only extra difficulty occurs in proving $k$-anti-uniformity for some $k\in \N$ that depends also on the degree of nilpotency of the group generated by  $T_1,\ldots, T_\ell$. In this case, the PET induction is somewhat more complicated, but can be handled by modifying the PET induction used in
 \mbox{\cite[Lemma 3.5]{FrHK11}} along the lines of the argument used to prove \mbox{\cite[Theorem 4.2]{W12}}.


 \subsection{Proof of Theorem~\ref{T:3}}
The inclusion $\overline{\mathcal{A}_{\ell}}\subset \overline{\mathcal{B}_\ell}$
follows from Proposition~\ref{P:nilkey}.
The inclusion $\overline{\mathcal{B}_{\ell}}\subset \overline{\mathcal{C}_\ell}$
is obvious.
The inclusion $\overline{\mathcal{C}_{\ell}}\subset \overline{\mathcal{A}_\ell}$
follows from Theorem~\ref{T:1}.


  \subsection{Proof of Theorems~\ref{T:4} and \ref{T:4'}}
The implication $\eqref{E:42}\Rightarrow \eqref{E:41}$ follows from Proposition~\ref{P:nilkey}.
(for Theorem~\ref{T:4'} in order to get property (i) for some fixed $\ell\in \N$ we use  property (ii) for $\ell!$).
The implication $\eqref{E:41}\Rightarrow \eqref{E:43}$ follows from Theorems~\ref{T:1} and \ref{T:1'}.
The implication $\eqref{E:43}\Rightarrow \eqref{E:42}$ is obvious.

 The same argument applies
for   the extensions mentioned after Theorem~\ref{T:4'} related to uniform and weighted Ces\`aro averages.

\end{document}